# FAMILIES OF NON-CONGRUENT NUMBERS

FRANZ LEMMERMEYER

ABSTRACT. In this article we study the Tate-Shafarevich groups corresponding to 2-isogenies of the curve $E_k : y^2 = x(x^2 - k^2)$ and construct infinitely many examples where these groups have odd 2-rank. Our main result is that among the curves $E_k$, where $k = pl \equiv 1 \bmod 8$ for primes $p$ and $l$, the curves with rank 0 have density $\geq \frac{1}{2}$.

## 1. INTRODUCTION

The elliptic curves $E_k : y^2 = x(x^2 - k^2)$ with $k \in \mathbb{Z}$ have been studied extensively, mainly because of the connection with the ancient problem of congruent numbers (see Guy [10] or Koblitz [14]). Many authors constructed families of non-congruent numbers by minimizing the Selmer groups attached to 2-isogenies of $E_k$ (see Feng [7, 8], Iskra [12], T. Ono [28], Serf [33], to name but the most recent contributors; actually results of this type go back to Genocchi [9] in the last century). Sharper results were obtained notably by J. Lagrange [16, 17] (see also Wada [36] and Nemenzo [25]), who found better bounds on the rank of $E_k$ by taking the 2-part of the Tate-Shafarevich groups into account. In this article, we will refine the criteria obtained by Lagrange and show that curves $E_k$, where $k = pl$ for primes $p \equiv l \equiv 1 \bmod 8$, very rarely have Tate-Shafarevich groups with trivial 2-part.

**Notation.** We recall the relevant notation from [18] (the standard reference for notions not explained here is Silverman [34]): elliptic curves $E$ with a rational point $T$ of order 2 as our curves $E_k$ come attached with a 2-isogeny $\phi : E \longrightarrow \widehat{E}$ (depending on the choice of $T$ if $E$ has three rational points of order 2). For $T = (0,0)$ we find the isogenous curve

$$\widehat{E}_k : y^2 = \begin{cases} x(x^2 + 4k^2) & \text{if } k \text{ is odd, and} \\ x(x^2 + k^2/4) & \text{if } k \text{ is even} \end{cases}$$

(the distinction is made in order to minimize the coefficients of the curve; we could just as well work with only $y^2 = x(x^2 + 4k^2)$ as both models are isomorphic). The dual isogeny $\widehat{E}_k \longrightarrow E_k$ will be denoted by $\psi$. If $k$ is fixed, we will suppress this index and write $E$ and $\widehat{E}$ for $E_k$ and $\widehat{E}_k$.

Consider the torsors

$$\mathcal{T}^{(\psi)}(b_1): \quad N^2 = b_1 M^4 + b_2 e^4, \quad b_1 b_2 = -k^2 \quad \text{and}$$
$$\mathcal{T}^{(\phi)}(b_1): \quad N^2 = b_1 M^4 + b_2 e^4, \quad b_1 b_2 = \begin{cases} 4k^2 & \text{if } k \text{ is odd,} \\ k^2/4 & \text{if } k \text{ is even.} \end{cases}$$

The main part of this article was written in 1999 while the author was at the MPI Bonn; he would like to thank everyone there for the hospitality and the stimulating environment, and the DFG for financial support during that time.

The Selmer group $S^{(\psi)}(\widehat{E}/\mathbb{Q})$ is defined as the subgroup of $\mathbb{Q}^\times/\mathbb{Q}^{\times 2}$ consisting of classes $b_1 \mathbb{Q}^{\times 2}$ such that $\mathcal{T}^{(\psi)}(b_1)$ has a nontrivial ($\neq (0,0,0)$) rational point in every completion $\mathbb{Q}_v$ of $\mathbb{Q}$; the subgroup of $S^{(\psi)}(\widehat{E}/\mathbb{Q})$ such that the torsors $\mathcal{T}^{(\psi)}(b_1)$ corresponding to $b_1 \mathbb{Q}^{\times 2}$ have a rational point will be denoted by $W(\widehat{E}/\mathbb{Q})$ (from now on, rational point will stand for non-trivial rational point; we may and do assume moreover that its coordinates are integral and primitive, that is, $(M, e) = 1$). Similarly we define $S^{(\phi)}(E/\mathbb{Q})$ and $W(E/\mathbb{Q})$. Finally, the Tate-Shafarevich groups are defined via the exact sequences

$$0 \longrightarrow W(E/\mathbb{Q}) \longrightarrow S^{(\phi)}(E/\mathbb{Q}) \longrightarrow \text{III}(E/\mathbb{Q})[\phi] \longrightarrow 0,$$

$$0 \longrightarrow W(\widehat{E}/\mathbb{Q}) \longrightarrow S^{(\psi)}(\widehat{E}/\mathbb{Q}) \longrightarrow \text{III}(\widehat{E}/\mathbb{Q})[\psi] \longrightarrow 0.$$

There exist various methods for constructing elements of order 2 in Tate-Shafarevich groups: one can perform a second 2-descent (cf. Birch and Swinnerton-Dyer [2], Razar [29], Lagrange [16, 17], Wada [36] and Nemenzo [24]), employ the Cassels pairing (see e.g. Aoki [1], Bölling [3], Cassels [4], and McGuinness [23]), compare the Selmer groups $S^{(\psi)}(\widehat{E}/\mathbb{Q})$ and $S^{(2)}(E/\mathbb{Q})$ as Kramer [15] (essentially, the methods mentioned so far are all equivalent to the classical second 2-descent), or use the method usually attributed to Lind [22] but actually going back (in a slightly different context) to Rédei [30] and Dirichlet [6] (I learned this technique from Stroeker & Top [35] and used it in [18] and [20]). In this paper, we continue to use this last method.

Our main results are the solvability criteria in Table 3 below; this will imply the lower bounds for density of rank-0 curves among the $E_{pl}$.

## 2. Preliminaries

In the calculations below we will have use quite a number of elementary results on quadratic reciprocity and genus theory. The following subsections recall what we will need.

**2.1. Some reciprocity laws.** In the following, $p$ and $l$ will denote primes $\equiv 1 \bmod 8$, and $\pi$ and $\lambda$ will denote primary primes in $\mathbb{Z}[i]$ with norms $p$ and $l$, respectively. A prime $\pi$ of norm $p \equiv 1 \bmod 8$ is called primary if $\pi$ is congruent to a square modulo 4. For $\pi \in \mathbb{Z}[i]$, $\Pi \in \mathbb{Z}[\sqrt{2}]$ and $\Pi^* \in \mathbb{Z}[\sqrt{-2}]$ we can always choose associates satisfying $\pi \equiv 1 \bmod 2 + 2i$, $\Pi \equiv 1 \bmod 2\sqrt{2}$ and $\Pi^* \equiv 1 \bmod 2\sqrt{-2}$, and these elements are primary.

We will need a few elementary results on quadratic residue symbols; as in [18], we let $(p/l)_4$ denote the biquadratic residue symbol for primes $l \equiv 1 \bmod 4$ such that $(p/l) = 1$, and we let $[\cdot/\cdot]$ denote the quadratic residue symbol in $\mathbb{Z}[i]$. We also note that, for primes $l = \lambda\overline{\lambda} \equiv 1 \bmod 8$, the relation $(1+i)^4 = -4$ implies that $[1+i/\lambda] = (-4/l)_8$ (this is the rational octic residue symbol). Moreover, $[\pi/\lambda] = (p/l)_4(l/p)_4$ for primes $p = \pi\overline{\pi}$ and $l = \lambda\overline{\lambda}$ such that $(p/l) = 1$ by Burde's rational reciprocity law. Finally, it is easy to check that $(\varepsilon_2/p) = [1+i/\pi] = (-4/p)_8$, where $\varepsilon_2 = 1 + \sqrt{2}$ (see [19]).

Now recall that primes $p \equiv 1 \bmod 8$ are norms from $\mathbb{Z}[\zeta_8]$, say $p = N\alpha$ for some $\alpha \equiv 1 \bmod (2 + 2\zeta)$, and in fact there exist primary elements $\pi \in \mathbb{Z}[i]$, $\Pi \in \mathbb{Z}[\sqrt{2}]$ and $\Pi^* \in \mathbb{Z}[\sqrt{-2}]$ with norm $p$. For primes $l \equiv 1 \bmod 8$, we define $\lambda$, $\Lambda$ and $\Lambda^*$ similarly. Unless explicitly stated otherwise, this notation is valid for the rest of this article.

The following result shows that solvability criteria involving the quadratic symbol $[\Pi^*/\Lambda^*]$ can be reduced to criteria involving only $[\Pi/\Lambda]$ and rational quartic residue symbols:

**Proposition 1.** *Let $p \equiv l \equiv 1 \mod 8$ be primes such that $(p/l) = +1$. Then*

$$\left[\frac{\Pi}{\Lambda}\right]\left[\frac{\Pi^*}{\Lambda^*}\right] = \left[\frac{\pi}{\lambda}\right] = \left(\frac{p}{l}\right)_4\left(\frac{l}{p}\right)_4,$$

*where the first three symbols $[\,\cdot\,/\,\cdot\,]$ denote the quadratic residue symbol in $\mathbb{Z}[\sqrt{2}\,]$, $\mathbb{Z}[\sqrt{-2}\,]$ and $\mathbb{Z}[i]$, respectively.*

*Proof.* We know that there exists an element $\alpha \in \mathbb{Z}[\zeta_8]$ such that $\Pi^* = \alpha_1\alpha_3$ (here $\alpha_j = \sigma_j(\alpha)$, where $\sigma_j$ is the automorphism that sends $\zeta_8$ to $\zeta_8^j$; in particular $\alpha_1 = \alpha$), $\pi = \alpha_1\alpha_5$ and $\Pi = \alpha_1\alpha_7$ (observe that such norms are necessarily totally positive). Defining $\beta$ accordingly we have $[\Pi/\Lambda] = (\alpha_1\alpha_7/\beta)$, where $(\,\cdot\,/\,\cdot\,)$ is the quadratic residue symbol in $\mathbb{Z}[\zeta_8]$. Similarly, we have $[\Pi^*/\Lambda^*] = (\alpha_1\alpha_3/\beta)$, hence $[\Pi/\Lambda][\Pi^*/\Lambda^*] = (\alpha_3\alpha_7/\beta)$. But this last symbol equals $[\overline{\pi}/\lambda]$, and since $(p/l) = +1$ this coincides with $[\pi/\lambda]$. This proves our claim by Burde's reciprocity law. $\square$

We also note that $[\Lambda/\Pi] = [\Pi/\Lambda]$ and $[\Lambda^*/\Pi^*] = [\Pi^*/\Lambda^*]$ by the quadratic reciprocity laws in $\mathbb{Z}[\sqrt{2}\,]$ and $\mathbb{Z}[\sqrt{-2}\,]$, respectively. Finally, if $K/k$ is an extension of number fields, if $[\,\cdot\,/\,\cdot\,]$ and $(\,\cdot\,/\,\cdot\,)$ denote the quadratic residue symbols in $K$ and $k$, respectively, and if $\mathfrak{a}$ is an ideal in $\mathcal{O}_k$ with odd norm, then $[\alpha/\mathfrak{a}] = (N_{K/k}\alpha/\mathfrak{a})$ right from the definition of residue symbols. Similarly, for ideals $\mathfrak{A}$ in $\mathcal{O}_K$ with relative norm $\mathfrak{a}$ and elements $\alpha \in k$ coprime to $\mathfrak{a}$, we have $[\alpha/\mathfrak{A}] = (\alpha/\mathfrak{a})$. For more on rational reciprocity laws, see [19, Chap. 5].

**2.2. The class groups of $\mathbb{Q}(\sqrt{\pm 2l}\,)$.** Let us begin by reviewing the basic results of Scholz as pertaining to the special case $k = \mathbb{Q}(\sqrt{2l}\,)$, where $l \equiv 1 \mod 4$ is prime. Let $\varepsilon$, $h$ and $h^+$ denote the fundamental unit, the class number and the class number in the strict sense of $k$. Moreover, define $(l/2)_4 = (-1/l)_8$ for primes $l \equiv 1 \mod 8$; then $(-4/l)_8 = (2/l)_4(l/2)_4$. The following proposition is the special case $p = 2$ of a more general result due to Scholz [31]:

**Proposition 2.** *With the notation as above, there are the following cases:*

- $(2/l) = -1$: *then $N\varepsilon = -1$ and $h \equiv h^+ \equiv 2 \mod 4$;*
- $(2/l) = +1$:
    (1) *if $(2/l)_4 = -(l/2)_4$, then $N\varepsilon = +1$, $h \equiv 2 \mod 4$, and $h^+ \equiv 4 \mod 8$.*
    (2) *if $(2/l)_4 = (l/2)_4 = -1$, then $N\varepsilon = -1$ and $h \equiv h^+ \equiv 4 \mod 8$;*
    (3) *if $(2/l)_4 = (l/2)_4 = +1$, then $4 \mid h$ and $8 \mid h^+$.*

Note that e.g. by the class number formula for strictly ambiguous ideals $C_{\text{am}} = 2^{t-1}/(E_F/NE_K)$ in quadratic extensions $K/F$ with $t$ ramified primes and unit groups $E_F$ and $E_K$, the prime ideal $\mathfrak{z}$ above 2 in $\mathbb{Q}(\sqrt{2}\,)$ is principal in the usual sense if and only if $N\varepsilon = +1$ for the fundamental unit $\varepsilon$ of $\mathbb{Q}(\sqrt{2}\,)$.

Let $\mathfrak{a} \stackrel{+}{\sim} \boxed{2}$ be short for "the ideal $\mathfrak{a}$ is equivalent in the strict sense to the square of some ideal", and define $\mathfrak{a} \stackrel{+}{\sim} \boxed{4}$ similarly.

If $d = d_1d_2$ is a product of two prime discriminants, then classical genus theory tells us that, for some ideal $\mathfrak{a}$ with norm $a$ (the existence of $\mathfrak{a}$ implies $(d/a) = +1$), we have $\mathfrak{a} \stackrel{+}{\sim} \boxed{2}$ if and only if $(d_1/a) = (d_2/a) = +1$.

**Lemma 3.** *Let $p \equiv l \equiv 1 \bmod 8$ be primes such that $(p/l) = +1$, and let $\mathfrak{p}$ denote the prime ideal above $p$ in $k = \mathbb{Q}(\sqrt{2l}\,)$. Then $\mathfrak{p} \stackrel{+}{\sim} \boxed{4} \iff [\Lambda/\Pi] = 1$.*

*Proof.* If $4 \mid h^+$, then the corresponding quartic cyclic unramified extension $K/k$ is given by $K = k(\sqrt{\Lambda}\,)$. A prime ideal $\mathfrak{p}$ of degree 1 will split completely in $K/k$ if and only if its ideal class is a fourth power in $\mathrm{Cl}^+(k)$; on the other hand, Kummer theory shows that it splits if and only if $\Lambda$ is a quadratic residue modulo any prime ideal above $\mathfrak{p}$ in $\mathbb{Q}(\sqrt{2}\,)$, that is, if and only if $[\Lambda/\Pi] = 1$.    □

**Lemma 4.** *Let $k = \mathbb{Q}(\sqrt{2l}\,)$ and assume that $(-4/l)_8 = -1$. Then the prime ideal $\mathfrak{z}$ above $2$ in $\mathcal{O}_k$ is principal in the strict sense if and only if $(2/l)_4 = -1$.*

*Proof.* First observe that our assumption implies by Proposition 2 that the fundamental unit of $k$ has positive norm, that $\mathfrak{z}$ is principal in the wide sense, and that $h^+ \equiv 4 \bmod 8$.

Assume that $\mathfrak{z}$ is principal in the strict sense. Then $X^2 - 2ly^2 = +2$ is solvable, hence so is $2x^2 - ly^2 = 1$ (we have put $X = 2x$). Now clearly $2 \nmid x$, hence $x^2 \equiv 1 \bmod 8$ and $2x^2 \equiv 2 \bmod 16$; on the other hand, $(2/y) = +1$, hence $y^2 \equiv 1 \bmod 16$. Together this implies that $l \equiv 1 \bmod 16$, that is, $(-1/l)_8 = +1$. Since $(-4/l)_8 = -1$ by assumption, this is equivalent to $(2/l)_4 = -1$.

Now assume that $\mathfrak{z}$ is not principal in the strict sense. Then $X^2 - 2ly^2 = -2$, and with $X = 2x$ we get $2x^2 - ly^2 = -1$. Now $(2/l)_4 = (x/l) = (l/x')$, where $x = 2^j x'$ with $x'$ odd, and $(l/x') = +1$ by reducing our equation modulo $x'$. Thus $(2/l)_4 = +1$.    □

## 3. The case $k = 2p$

We will now investigate which torsors of $E_{2p}$ do not have rational points although they are everywhere locally solvable. These curves were already studied by Lagrange [17] using second 2-descents and by Kings [13] using the Cassels pairing on $\mathrm{III}(E/\mathbb{Q})$. The curves $E_{2p}$ are the simplest examples where $\mathrm{III}(E/\mathbb{Q})[\phi]$ and $\mathrm{III}(\widehat{E}/\mathbb{Q})[\psi]$ may have odd dimension:

**Theorem 5.** *Let $p \equiv 1 \bmod 8$ be a prime and consider the elliptic curve $E : y^2 = x(x^2 - 4p^2)$. Then the Selmer groups are given by*

$$S^{(\psi)}(\widehat{E}/\mathbb{Q}) = \langle -1, 2, p \rangle, \quad S^{(\phi)}(E/\mathbb{Q}) = \langle p \rangle,$$

*and if $p \equiv 9 \bmod 16$, then we have $\mathrm{III}(\widehat{E}/\mathbb{Q})[\psi] = \langle p \rangle$ and $\mathrm{III}(E/\mathbb{Q})[\phi] = \langle p \rangle$. Here $\langle x, \ldots, z \rangle$ denotes the subgroup of $\mathbb{Q}^\times/\mathbb{Q}^{\times 2}$ generated by $x, \ldots, z$. Moreover, $\mathrm{III}(E/\mathbb{Q})[2] \simeq \mathrm{III}(\widehat{E}/\mathbb{Q})[2] \simeq (\mathbb{Z}/2\mathbb{Z})^2$.*

*Proof.* We leave the proofs that $\mathrm{III}(\widehat{E}/\mathbb{Q})[\psi]$ and $\mathrm{III}(E/\mathbb{Q})[\phi]$ both have order 2 as an exercise to the reader (they are much simpler than the proofs in the sections below). The claims $\mathrm{III}(E/\mathbb{Q})[2] \simeq \mathrm{III}(\widehat{E}/\mathbb{Q})[2] \simeq (\mathbb{Z}/2\mathbb{Z})^2$ follow from the exact sequences (see diagram (3.9) in Razar [29, p. 139]; Feng [7, 8] erroneously claims that $C = \widehat{C} = 0$)

$$0 \longrightarrow \mathrm{III}(\widehat{E}/\mathbb{Q})[\psi] \longrightarrow \mathrm{III}(\widehat{E}/\mathbb{Q})[2] \longrightarrow \mathrm{III}(E/\mathbb{Q})[\phi] \longrightarrow C \longrightarrow 0$$
$$0 \longrightarrow \mathrm{III}(E/\mathbb{Q})[\phi] \longrightarrow \mathrm{III}(E/\mathbb{Q})[2] \longrightarrow \mathrm{III}(\widehat{E}/\mathbb{Q})[\psi] \longrightarrow \widehat{C} \longrightarrow 0$$

where $C$ and $\widehat{C}$ are finite groups of square order by results of Cassels [?]. Since they are quotients of groups of order 2, it follows that $C = \widehat{C} = 0$, and this implies our claim. □

## 4. The case $k = pl \equiv 1 \bmod 8$

The simplest cases are those where $p \equiv l \equiv 3, 5, 7 \bmod 8$; they were already discussed by Lagrange [16]:

| $p \bmod 8$ | $l \bmod 8$ | $(p/l)$ | $S^{(\psi)}(\widehat{E}/\mathbb{Q})$ | $S^{(\phi)}(E/\mathbb{Q})$ |
|---|---|---|---|---|
| 1 | 1 | $+1$ | $\langle -1, p, l \rangle$ | $\langle 2, p, l \rangle$ |
|   |   | $-1$ | $\langle -1, pl \rangle$ | $\langle 2, pl \rangle$ |
| 5 | 5 | $+1$ | $\langle -1, pl \rangle$ | $\langle p, l \rangle$ |
|   |   | $-1$ | $\langle -1, pl \rangle$ | $\langle 2p, 2l \rangle$ |
| 3 | 3 |   | $\langle -1, pl \rangle$ | 1 |
| 7 | 7 |   | $\langle -1, p, l \rangle$ | $\langle 2 \rangle$ |

He also found necessary criteria for the solvability of certain torsors. Here are the results, reformulated using our notation:

**Proposition 6.** *Let $p$ and $l$ be distinct primes such that $p \equiv l \equiv 3, 5, 7 \bmod 8$. If the torsors in the table below have a rational point, then the conditions in the last column must be satisfied:*

| $p \bmod 8$ | $l \bmod 8$ | $(p/l)$ | torsors | conditions |
|---|---|---|---|---|
| 5 | 5 | $+1$ | $\mathcal{T}^{(\phi)}(p), \mathcal{T}^{(\phi)}(l), \mathcal{T}^{(\phi)}(pl)$ | $(p/l)_4 = (l/p)_4$ |
| 5 | 5 | $-1$ | $\mathcal{T}^{(\phi)}(2p), \mathcal{T}^{(\phi)}(2l), \mathcal{T}^{(\phi)}(pl)$ | $[(1+i)\pi/\lambda] = +1$ |
| 7 | 7 | $+1$ | $\mathcal{T}^{(\phi)}(2), \mathcal{T}^{(\psi)}(p), \mathcal{T}^{(\psi)}(l)$ | $[\Lambda/\Pi] = +1$ |

*In the last row, $\Lambda \in \mathbb{Z}[\sqrt{2}\,]$ is a primary element with norm $-l$, and $\Pi \in \mathbb{Z}[\sqrt{2}\,]$ has norm $\pm p$. Observe that $[\Lambda/\Pi]$ is well defined since $[\Lambda/\Pi][\overline{\Lambda}/\Pi] = [-l/\Pi] = (-l/p) = (p/l) = +1$.*

The proofs for $p \equiv l \equiv 5 \bmod 8$ are straight forward and left as an exercise to the reader. Here we give some details for the case $p \equiv l \equiv 7 \bmod 8$: Consider the torsor $\mathcal{T}(p) : pn^2 = M^4 - l^2 e^4$. Reduction modulo $l$ shows immediately that either 1) $l \nmid M$ and $(p/l) = +1$, or 2) $l \mid M$ and $(p/l) = -1$. Moreover, either A) $2 \nmid Me$ and $2 \mid n$, or B) $2 \nmid ne$ and $2 \mid M$. As in the case $p \equiv l \equiv 1 \bmod 8$, we get four equations per case:

| case | I | II | III | IV |
|---|---|---|---|---|
| 1A) | $M^2 + le^2 = 2pa^2$ | $M^2 - le^2 = 2b^2$ | $pa^2 + b^2 = M^2$ | $pa^2 - b^2 = le^2$ |
| 1B) | $M^2 + le^2 = pa^2$ | $M^2 - le^2 = b^2$ | $pa^2 + b^2 = 2M^2$ | $pa^2 - b^2 = 2le^2$ |
| 2A) | $lm^2 + e^2 = 2a^2$ | $lm^2 - e^2 = 2pb^2$ | $a^2 - pb^2 = e^2$ | $a^2 - pb^2 = lm^2$ |
| 2B) | $lm^2 + e^2 = a^2$ | $lm^2 - e^2 = pb^2$ | $a^2 - pb^2 = 2e^2$ | $a^2 - pb^2 = 2lm^2$ |

Now we distinguish these four cases:

1A) Writing II) and III) in the form $M^2 - 2b^2 = le^2$ and $M^2 - b^2 = pa^2$ we find that $[\lambda/\Pi] = [M + b\sqrt{2}/\Pi]$, where $\lambda \in \mathbb{Z}[\sqrt{2}]$ is the element of norm $l$ that divides $M + b\sqrt{2}$. We would like to use the congruence $b \equiv \pm M \bmod p$ coming from second equation and conclude that $[\lambda/\Pi] = [1 \pm \sqrt{2}/\Pi]$, but unfortunately the last symbol depends on the choice of the sign. We therefore have to work a little harder.

First observe that, for $M > 0$, we have $(M/p) = (-p/M) = +1$ from the second equation. Now we factor $pa^2 = (M - b)(M + b)$ and consider the following two cases:

a) $M - b = 2pr^2$, $M + b = 2s^2$ (the negative signs cannot hold here: otherwise we would get $M = -s^2 - rp^2$ contradicting our assumption that $M > 0$); then $b = s^2 - pr^2 \equiv 1 \bmod 4$, and we get $[M + b\sqrt{2}/\Pi] = (M/p)[1 + \sqrt{2}/\Pi] = [1 + \sqrt{2}/\Pi]$.

b) $M - b = 2r^2$, $M + b = 2ps^2$; then $b = ps^2 - r^2 \equiv 3 \bmod 4$, and now $[M + b\sqrt{2}/\Pi] = [1 - \sqrt{2}/\Pi]$.

Thus $[\lambda\varepsilon/\Pi] = +1$, where $\varepsilon = 1 + (-1/b)\sqrt{2}$. Now it is easy to check that $\lambda\varepsilon$ is primary: in fact, $r + s\sqrt{2}$ with $2 \mid s$ is primary if and only if $r + s \equiv 1 \bmod 4$, and since $\lambda\varepsilon$ is primary if and only if $(M + b\sqrt{2})\varepsilon$ is, we find

$$(M + b\sqrt{2})\varepsilon = \begin{cases} M + 2b + (M + b)\sqrt{2} & \text{if } b \equiv 1 \bmod 4, \\ M - 2b + (b - M)\sqrt{2} & \text{if } b \equiv 3 \bmod 4, \end{cases}$$

and $2M + 3b \equiv 2 - b \equiv 1 \bmod 4$ in the first and $-b \equiv 1 \bmod 4$ in the second case. Thus in this case $[\Lambda/\Pi] = +1$, where $\Lambda = \lambda\varepsilon$ is primary with norm $-l$.

1B) Here $b^2 - 2M^2 = -pa^2$ and $b^2 - M^2 = -le^2$. Again, choosing $M > 0$ guarantees $(M/l) = (-l/M) = +1$. Next, $M - b = r^2$ and $M + b = ls^2$ imply $2b = ls^2 - r^2$ and $b \equiv 1 \bmod 4$, while $M - b = lr^2$ and $M + b = s^2$ give $2b = s^2 - lr^2$ and $b \equiv 3 \bmod 4$. Thus we get $[b + M\sqrt{2}/\Lambda] = [-1 + \sqrt{2}/\Lambda]$ if $b \equiv 1 \bmod 4$ and $[b + M\sqrt{2}/\Lambda] = [1 + \sqrt{2}/\Lambda]$ if $b \equiv 3 \bmod 4$. Putting $\varepsilon = -(-1/b) + \sqrt{2}$, it is easy to check that $(b + M\sqrt{2})\varepsilon$ is totally positive. Now Hasse [11] has shown that we have the reciprocity law $[\alpha/\beta] = [\beta/\alpha]$ in an arbitrary algebraic number field if the conductors of $\alpha$ and $\beta$ are coprime. Since $(b + M\sqrt{2})\varepsilon \gg 0$, the gcd of the conductors of $(b + M\sqrt{2})\varepsilon$ and $\Lambda$ do not contain infinite primes, and since $\Lambda$ is primary, the gcd does not contain primes above 2. But then $(b + M\sqrt{2}, \Lambda) = (1)$ gurantees that the conductors are indeed coprime, and the reciprocity law gives $1 = [(b + M\sqrt{2})\varepsilon/\Lambda] = [\Lambda/b + M\sqrt{2}]$. Since $(b + M\sqrt{2}) = (\Pi\alpha^2)$ by unique factorization, we conclude that $[\Lambda/b + M\sqrt{2}] = [\Lambda/\Pi]$.

2A) Equations I and III correspond to III and II in case 1B) with the roles of $p$ and $l$ switched.

2B) Again, this reduces to case 1A).

We have proved:

**Proposition 7.** *Let $d_i$ denote the density of rank $0$ curves among the $E_{pl}$, where $p \equiv l \equiv i \bmod 8$ are primes. Then we have $d_3 = 1$, $d_5 \geq \frac{1}{2}$ and $d_7 \geq \frac{1}{2}$.*

The main result of this note is that we also have $d_1 \geq \frac{1}{2}$ (this is much stronger than the result obtained by Lagrange [17]). Although numerical computations seem to suggest that $d_i = 1$, it seems that the bounds derived in this article cannot be improved using our methods.

From now on, we will assume that $p$ and $q$ are both primes $\equiv 1 \mod 8$.

4.1. **The case** $(p/l) = -1$. Let $k = pl$ be a product of primes $p \equiv l \equiv 1 \mod 8$ with $(p/l) = -1$. Then (see [17])

$$S^{(\psi)}(\widehat{E}/\mathbb{Q}) = \langle -1, pl \rangle = W(\widehat{E}/\mathbb{Q}), \quad S^{(\phi)}(E/\mathbb{Q}) = \langle 2, pl \rangle.$$

In particular, $\text{III}(\widehat{E}/\mathbb{Q})[\psi] = 0$, so we only have to discuss the $\phi$-part of $\text{III}(E/\mathbb{Q})$.

**Proposition 8.** *If $k = pl$ is a product of primes $p \equiv l \equiv 1 \mod 8$ with $(p/l) = -1$, then $\text{III}(E/\mathbb{Q})[\phi] = \langle 2, pl \rangle$ whenever $(-4/p)_8(-4/l)_8 = -1$. If this condition holds, we have $\#\text{III}(E/\mathbb{Q})[2] = 4$.*

*Proof.* Consider $\mathcal{T}^{(\phi)}(2) : N^2 = 2M^4 + 2p^2l^2e^4$.

- Assume first that $(M, pl) = 1$; then $N = 2n$ gives $2n^2 = M^4 + p^2l^2e^4$. Now $M^2 + ple^2i \equiv 1 + i \mod 8$ and unique factorization in $\mathbb{Z}[i]$ shows that $M^2 + ple^2i = (1+i)\nu^2$. Write $p = \pi\overline{\pi}$ for primes $\pi, \overline{\pi} \equiv 1 \mod 2 + 2i$; reducing modulo $\pi$ gives $[1 + i/\pi] = +1$, that is, $(-4/p)_8 = +1$, and similarly $(-4/l)_8 = +1$.
- If $(M, pl) = p$, put $M = mp$ and $N = 2pn$; then we get $2n^2 = (pm^2 + le^2i)(pm^2 + le^2i)$, and again $pm^2 + le^2i = (1+i)\nu^2$. Reducing modulo $\pi$ gives $[1 + i/\pi] = [l/\pi] = (l/p) = -1$, hence $(-4/p)_8 = (-4/l)_8 = -1$.
- The cases $(M, pl) = l$ and $(M, pl) = pl$ are treated similarly.

Next take $\mathcal{T}^{(\phi)}(pl) : N^2 = plM^4 + 4ple^4$. With $N = pln$ this gives $pln^2 = M^4 + 4e^4$; since we may switch the roles of $M$ and $e$ we may assume that $M$ is odd and $e$ is even. Reducing modulo $p$ and $l$ shows that $(-4/pl)_8 = (Me/p)$. Write $e = 2^j e'$ with $e'$ odd: then $(e/p) = (e'/p) = (p/e') = 1$ and $(M/p) = (p/M) = 1$. Thus $(-4/pl)_8 = 1$.

Finally look at $\mathcal{T}^{(\phi)}(2pl) : 2pln^2 = M^4 + e^4$. As above, $M^2 + ie^2 = (1+i)\pi\lambda\nu^2$; adding this equation to its conjugate gives $2M^2 = (1+i)\pi\lambda\nu^2 + (1-i)\overline{\pi}\overline{\lambda}\overline{\nu}^2$. Reducing modulo $\overline{\pi}$ gives $1 = (2/p) = [1 + i/\overline{\pi}][\pi/\overline{\pi}][\lambda/\overline{\pi}]$. Now $[\pi/\overline{\pi}] = 1$ and $[\lambda/\overline{\pi}] = [\lambda/\pi]$, hence $(-4/p)_8 = [\pi/\lambda]$. Similarly, $(-4/l)_8 = [\pi/\lambda]$, and the claim follows. Note that $[\pi/\lambda]$ depends on the choice of $\pi$ and $\lambda$. □

From Prop. 8 we get by a standard application of Chebotarev's density theorem the following

**Corollary 9.** *The curves of rank $0$ among $E_{pl}$, where $p \equiv l \equiv 1 \mod 8$ are primes such that $(p/l) = -1$, have density at least $\frac{1}{2}$.*

4.2. **The case** $(p/l) = +1$. Let $k = pl$ be a product of primes $p \equiv l \equiv 1 \mod 8$ with $(p/l) = +1$. Then (see [17])

$$S^{(\psi)}(\widehat{E}/\mathbb{Q}) = \langle -1, p, l \rangle, \quad S^{(\phi)}(E/\mathbb{Q}) = \langle 2, p, l \rangle.$$

Moreover $\langle -1, pl \rangle \subseteq W(\widehat{E}/\mathbb{Q})$. As above, we will now compute nontrivial elements of $\text{III}(E/\mathbb{Q})[\phi]$ and $\text{III}(\widehat{E}/\mathbb{Q})[\psi]$.

$$\boxed{\text{The } \psi\text{-part}}$$

First we observe that $W(\widehat{E}/\mathbb{Q})$ always contains $\langle -1, pl \rangle$. Thus

$$\text{either} \quad W(\widehat{E}/\mathbb{Q}) = \langle -1, p, l \rangle \text{ and } \mathbf{III}(\widehat{E}/\mathbb{Q})[\psi] = 0,$$
$$\text{or} \quad W(\widehat{E}/\mathbb{Q}) = \langle -1, pl \rangle \text{ and } \mathbf{III}(\widehat{E}/\mathbb{Q})[\psi] = \langle p \rangle,$$

where $\langle p \rangle$ represents the class of $p\mathbb{Q}^{\times 2}$ (which is the same as the class of $l\mathbb{Q}^{\times 2}$ in view of $pl\mathbb{Q}^{\times 2} \in W(\widehat{E}/\mathbb{Q})$) in $\mathbf{III}(\widehat{E}/\mathbb{Q})[\psi]$.

It is therefore sufficient to consider the torsor $\mathcal{T}^{(\psi)}(p) : N^2 = pM^4 - pl^2 e^4$. Here the right hand side factors over $\mathbb{Q}$ as $N^2 = p(M^2 - le^2)(M^2 + le^2)$. We have the following possibilities concerning divisibility by 2: $\begin{cases} 1) \ 2 \mid e, 2 \nmid MN \\ 2) \ 2 \mid N, 2 \nmid Me, \end{cases}$ by $l$: $\begin{cases} A) \ l \nmid MN \\ B) \ l \mid M, l \mid N, \end{cases}$ and by $p$: $\begin{cases} a) \ p \mid (M^2 + le^2) \\ b) \ p \mid (M^2 - le^2). \end{cases}$

Thus we have to consider eight different cases. We claim

**Proposition 10.** *Let $E$ be the elliptic curve defined by $y^2 = x(x^2 - k^2)$, where $k = pl$ and where $p \equiv l \equiv 1 \bmod 8$ are primes such that $(p/l) = 1$. If the torsor*

$$\mathcal{T}^{(\psi)}(p) : N^2 = pM^4 - pl^2 e^4, \tag{1}$$

*has a rational solution, then the conditions in Table 1 hold according to the case we are in.*

| case | conditions($*$) |
|---|---|
| 1Aa) | $[\Pi/\Lambda] = (l/p)_4 = (-4/p)_8 = 1$ |
| 1Ab) | $[\Pi/\Lambda] = (p/l)_4 = (l/p)_4(-4/p)_8 = 1$ |
| 1Ba) | $[\Pi/\Lambda] = (-4/pl)_8, (l/p)_4 = (p/l)_4(-4/p)_8 = 1$ |
| 1Bb) | $[\Pi/\Lambda] = (-4/l)_8, (p/l)_4 = (-4/p)_8 = 1$ |
| 2Aa) | $[\Pi/\Lambda] = (-4/p)_8, (l/p)_4 = (-4/l)_8 = 1$ |
| 2Ab) | $[\Pi/\Lambda] = (l/p)_4 = (p/l)_4(-4/l)_8 = 1$ |
| 2Ba) | $[\Pi/\Lambda] = (-4/pl)_8, (p/l)_4 = (l/p)_4(-4/l)_8 = 1$ |
| 2Bb) | $[\Pi/\Lambda] = (p/l)_4 = (-4/l)_8 = 1$ |

TABLE 1. If $\mathcal{T}^{(\psi)}(p)$ has a rational point, then the conditions ($*$) hold.

If we are in case 1A), then putting $N = pn$ in (1) gives $pn^2 = M^4 - l^2 e^4 = (M^2 - le^2)(M^2 + le^2)$. In case 1Aa), these two factors are coprime, hence $M^2 + le^2 = pa^2$ (I) and $M^2 - le^2 = b^2$ (II), where $ab = n$. By adding and subtracting (I) and (II) we get $2M^2 = b^2 + pa^2$ (III) and $2le^2 = pa^2 - b^2$ (IV). In a similar way we find the following table displaying the four equations (I)–(IV) whose solvability follow from the existence of a rational point on (1):

| case | I | II | III | IV |
|---|---|---|---|---|
| $1Aa)$ | $M^2 + le^2 = pa^2$ | $M^2 - le^2 = b^2$ | $2M^2 = b^2 + pa^2$ | $2le^2 = pa^2 - b^2$ |
| $1Ab)$ | $M^2 + le^2 = a^2$ | $M^2 - le^2 = pb^2$ | $2M^2 = a^2 + pb^2$ | $2le^2 = a^2 - pb^2$ |
| $1Ba)$ | $lm^2 + e^2 = pa^2$ | $lm^2 - e^2 = b^2$ | $2lm^2 = b^2 + pa^2$ | $2e^2 = pa^2 - b^2$ |
| $1Bb)$ | $lm^2 + e^2 = a^2$ | $lm^2 - e^2 = pb^2$ | $2lm^2 = a^2 + pb^2$ | $2e^2 = a^2 - pb^2$ |
| $2Aa)$ | $M^2 + le^2 = 2pa^2$ | $M^2 - le^2 = 2b^2$ | $M^2 = b^2 + pa^2$ | $le^2 = pa^2 - b^2$ |
| $2Ab)$ | $M^2 + le^2 = 2a^2$ | $M^2 - le^2 = 2pb^2$ | $M^2 = a^2 + pb^2$ | $le^2 = a^2 - pb^2$ |
| $2Ba)$ | $lm^2 + e^2 = 2pa^2$ | $lm^2 - e^2 = 2b^2$ | $lm^2 = b^2 + pa^2$ | $e^2 = pa^2 - b^2$ |
| $2Bb)$ | $lm^2 + e^2 = 2a^2$ | $lm^2 - e^2 = 2pb^2$ | $lm^2 = a^2 + pb^2$ | $e^2 = a^2 - pb^2$ |

In order to save some work we prove a general result that may be applied to each of these cases:

**Proposition 11.** *Let $A, B, C, D \in \mathbb{N}$ be pairwise coprime integers, each a product of primes $\equiv 1 \bmod 4$, and assume that these primes are quadratic residues of each other. If there are $x, y, v, w \in \mathbb{N}$ such that*

$$Ax^2 + By^2 = Cv^2, \tag{2}$$

$$Ax^2 - By^2 = Dw^2, \tag{3}$$

*then $C \equiv D \bmod 8$, and $A, B, C$ and $D$ satisfy the relations*

$$\Big(\frac{AB}{C}\Big)_4 \Big(\frac{AD}{B}\Big)_4 \Big(\frac{BD}{A}\Big)_4 = 1 \tag{4}$$

*and*

$$(-1)^{\frac{C-D}{8}} \Big(\frac{2}{CD}\Big)_4 \Big(\frac{BC}{D}\Big)_4 \Big(\frac{BD}{C}\Big)_4 \Big(\frac{CD}{A}\Big)_4 = 1. \tag{5}$$

*Proof.* Assume that we have a congruence $Ar^2 = Bs^2 \bmod D$ with $(r, D) = (s, D) = 1$, and assume moreover that $(AB/p) = +1$ for all $p \mid D$. Then for each such $p$ we have $Ar^2 = Bs^2 \bmod p$, and raising this congruence to the $\frac{p-1}{4}$-th power we find that $(A/p)_4 (r/p) = (B/p)_4 (s/p)$; multiplying these relations together shows that $(AB/D)_4 = (rs/D)$. We will use this type of reasoning without comment below.

We may (and will) assume that $(x, y) = 1$. From $2y^2 \equiv 2By^2 = Cv^2 - Dw^2 \equiv v^2 - w^2 \bmod 4$ we then deduce that $2 \mid y$ and $2 \nmid xvw$.

Reducing (2) modulo $C$ gives $(-AB/C)_4 = (xy/C)$. Writing $y = 2^j y'$ for some odd $y'$ gives $(y/C) = (2/C)^j (y'/C) = (2/C)^j (C/y')$. Reducing (2) modulo $y'$ we see $(C/y') = (A/y')$. Similarly, we get $(x/C) = (C/x) = (B/x) = (x/B)$ from (2), and $(x/B) = (AD/B)_4 (w/B)$. Since $(w/B) = (B/w) = (A/w) = (w/A) = (-BD/A)_4 (y/A) = (-BD/A)_4 (2/A)^j (A/y')$, collecting our results gives $(-AB/C)_4 = (AD/B)_4 (-BD/A)_4 (2/AC)^j$. Next, $(-1/A)_4 = (2/A)$ and $(-1/C) = (2/C)$, hence the relations becomes $(AB/C)_4 (AD/B)_4 (BD/A)_4 = (2/AC)^{j+1}$.

Now there are two cases: if $j = 1$, then $A \equiv C + 4 \bmod 8$, hence $(2/AC) = -1$, but $(2/AC)^{j+1} = 1$; if $j \geq 2$, then $A \equiv C \bmod 8$, hence $(2/AC) = 1$. In both cases, we arrive at the desired relation.

By adding and subtracting (2) and (3), we get

(6) $$2Ax^2 = Cv^2 + Dw^2,$$

(7) $$2By^2 = Cv^2 - Dw^2.$$

From (7) and the fact that $y$ is even we deduce that $C \equiv D \bmod 8$.

Reducing (7) modulo $D$ yields $(2BC/D)_4 = (vy/D)$. From (7) we deduce that $(y/D) = (2/D)^j(y'/D) = (2/D)^j(D/y') = (2/D)^j(C/y')$ and $(v/D) = (D/v) = (2A/v)$, so $(2BC/D)_4 = (2A/v)(2/D)^j(C/y')$. Similarly, $(-2BD/C)_4 = (wy/C)$, $(y/C) = (2/C)^j(C/y')$ and $(w/C) = (C/w) = (2A/w)$. Combining these results yields $(2BC/D)_4(-2BD/C)_4 = (2A/vw)(2/CD)^j$. Since $C \equiv D \bmod 8$, we have $(2/CD) = 1$, and using $(-1/C)_4 = (2/C)$ we conclude that

$$\Big(\frac{2BC}{D}\Big)_4 \Big(\frac{2BD}{C}\Big)_4 = \Big(\frac{2A}{vw}\Big)\Big(\frac{2}{C}\Big).$$

Next, $(A/w) = (w/A) = (-BD/A)_4(y/A)$ and $(A/v) = (v/A) = (BC/A)_4(y/A)$, thus $(A/vw) = (-BD/A)_4(BC/A)_4 = (2/A)(CD/A)_4$ since $(B/A) = +1$. This gives us

$$\Big(\frac{2}{CD}\Big)_4 \Big(\frac{BC}{D}\Big)_4 \Big(\frac{BD}{C}\Big)_4 \Big(\frac{CD}{A}\Big)_4 = \Big(\frac{2}{vw}\Big)\Big(\frac{2}{C}\Big).$$

If $j = 1$, then $Cv^2 \equiv Dw^2 + 8 \bmod 16$, hence $C \equiv D + 8 \bmod 16$ if and only if $(2/v) = (2/w)$, or $(2/vw) = -(-1)^{(C-D)/8}$. Moreover, $2Ax^2 \equiv 2Cv^2 + 8 \bmod 16$ implies $(2/AC) = -1$, so we get $(2/vw)(2/AC) = (-1)^{(C-D)/8}$.

If $j \geq 2$, then $Cv^2 \equiv Dw^2 \bmod 16$, and this shows that $C \equiv D \bmod 16$ if and only if $(2/v) = (2/w)$, hence $(2/vw) = (-1)^{(C-D)/8}$. Moreover, (6) implies that $A \equiv C \bmod 8$, hence $(2/AC) = +1$, and again $(2/vw)(2/AC) = (-1)^{(C-D)/8}$. $\square$

In order to apply this result we have to identify the coefficients $A, B, C$ and $D$. We find

| case | (1) | (2) | A | B | C | D | case | (1) | (2) | A | B | C | D |
|------|-----|-----|---|---|---|---|------|-----|-----|---|---|---|---|
| 1Aa) | I | II | 1 | $l$ | $p$ | 1 | 2Aa) | III | IV | $p$ | 1 | 1 | $l$ |
| 1Ab) | I | II | 1 | $l$ | 1 | $p$ | 2Ab) | III | IV | 1 | $p$ | 1 | $l$ |
| 1Ba) | I | II | $l$ | 1 | $p$ | 1 | 2Ba) | III | IV | $p$ | 1 | $l$ | 1 |
| 1Bb) | I | II | $l$ | 1 | 1 | $p$ | 2Bb) | III | IV | 1 | $p$ | $l$ | 1 |

This takes care of all the conditions not involving $[\Pi/\Lambda]$. For completing the proof we need the following

**Lemma 12.** *Let $P \equiv L \equiv 1 \bmod 8$ be primes such that $(P/L) = +1$. Let $\Pi, \Lambda \in \mathbb{Z}[\sqrt{2}\,]$ be primary elements of norm $P$ and $L$, respectively. If there exist integers $x, y, z, w \in \mathbb{N}$ such that*

$$x^2 - 2y^2 = -Pz^2, \quad \text{and} \quad x^2 - y^2 = \epsilon L w^2$$

*for some $\epsilon = \pm 1$, then $[\Pi/\Lambda] = +1$.*

*Proof.* Unique factorization gives $x + y\sqrt{2} = \varepsilon_2 \Pi \alpha^2$, where $\varepsilon_2$ is a fundamental unit of $\mathbb{Z}[\sqrt{2}\,]$ and where $N\alpha = z$. Thus $[\Pi/\Lambda] = [\varepsilon_2/\Lambda][x + y\sqrt{2}\Lambda]$. Now $y \equiv \pm x \bmod \Lambda$ from the second equation, hence $[x + y\sqrt{2}/\Lambda] = [x/\Lambda][1 + \sqrt{2}/\Lambda]$. But $[1 + \sqrt{2}/\Lambda] =$ $[\varepsilon_2/\Lambda]$ since the expression $[\pm 1 \pm \sqrt{2}/\Lambda]$ does not depend on the choice of signs, and we get $[\Pi/\Lambda] = [x/\Lambda] = (x/L)$. If $\epsilon = +1$, then $(x/L) = (y/L) = (L/y) = +1$, and if $\epsilon = -1$, then $(x/L) = (L/x) = +1$. This proves our claim. □

Lemma 12 takes care of four out of our eight cases:

| case | $x$ | $y$ | $z$ | $w$ | $P$ | $L$ | $\epsilon$ |
|---|---|---|---|---|---|---|---|
| 1Aa) | $b$ | $M$ | $a$ | $e$ | $p$ | $l$ | $-1$ |
| 1Ab) | $a$ | $M$ | $b$ | $e$ | $p$ | $l$ | $+1$ |
| 2Ab) | $M$ | $a$ | $e$ | $b$ | $l$ | $p$ | $+1$ |
| 2Bb) | $e$ | $a$ | $m$ | $b$ | $l$ | $p$ | $-1$ |

For the remaining four cases, the role of Lemma 12 is taken over by

**Lemma 13.** *Let $P \equiv L \equiv 1 \bmod 8$ be primes such that $(P/L) = +1$. Let $\Pi, \Lambda \in \mathbb{Z}[\sqrt{2}]$ be primary elements of norm $P$ and $L$, respectively. If there exist integers $x, y, z, w \in \mathbb{N}$ such that*

$$x^2 + 2\epsilon y^2 = Pz^2, \quad \text{and} \quad x^2 + \epsilon y^2 = Lw^2$$

*for some $\epsilon = \pm 1$, then*

$$\left[\frac{\Pi}{\Lambda}\right] = \begin{cases} \left(\frac{-4}{L}\right)_8 & \text{if } \epsilon = -1, \\ \left(\frac{P}{L}\right)_4 \left(\frac{L}{P}\right)_4 \left(\frac{-4}{L}\right)_8 & \text{if } \epsilon = +1. \end{cases}$$

*Proof.* Let $\pi, \lambda \in \mathbb{Z}[\sqrt{2\epsilon}]$ be primary elements of norm $P$ and $L$, respectively. Then from $\pi\alpha^2 = x + y\sqrt{2\epsilon}$ we get $[\pi/\lambda] = [x + y\sqrt{2\epsilon}/\lambda]$. The second equation gives $x \equiv \pm y\sqrt{\epsilon} \bmod \mathfrak{l}$, where $\mathfrak{l}$ denotes a prime ideal above $l$ in $\mathbb{Q}(\zeta_8)$. Letting $\{\cdot/\cdot\}$ denote the quadratic residue symbol in $\mathbb{Z}[\zeta_8]$, we find $[x+y\sqrt{2\epsilon}/\lambda] = \{x+y\sqrt{2\epsilon}/\mathfrak{l}\} = \{x \pm x\sqrt{2}/\mathfrak{l}\} = (x/L)[1 \pm \sqrt{2}/\Lambda]$. Now if $\epsilon = 1$ then $(x/L) = (L/x) = +1$, whereas if $\epsilon = -1$ then $(x/L) = (y/L) = (y'/L) = (L/y') = +1$. Thus $[\pi/\lambda] = [1+\sqrt{2}/\Lambda] = (-4/L)_8$. If $\epsilon = -1$, then $\pi = \Pi$ and $\lambda = \Lambda$, but if $\epsilon = +1$ then $\pi = \Pi^*$ and $\lambda = \Lambda^*$, $\Pi^*, \Lambda^* \in \mathbb{Z}[\sqrt{-2}]$ are primary elements of norm $p$ and $l$, respectively. Thus $[\Pi/\Lambda] = [\Pi^*/\Lambda^*](P/L)_4(L/P)_4 = (P/L)_4(L/P)_4(-4/L)_8$. □

Lemma 13 covers the remaining four cases:

| case | $x$ | $y$ | $z$ | $w$ | $P$ | $L$ | $\epsilon$ | resulting condition |
|---|---|---|---|---|---|---|---|---|
| 1Ba) | $b$ | $e$ | $a$ | $m$ | $p$ | $l$ | $+1$ | $[\Pi/\Lambda] = (-4/pl)_8$ |
| 1Bb) | $a$ | $e$ | $m$ | $b$ | $p$ | $l$ | $-1$ | $[\Pi/\Lambda] = (-4/l)_8$ |
| 2Aa) | $M$ | $b$ | $a$ | $e$ | $l$ | $p$ | $-1$ | $[\Pi/\Lambda] = (-4/p)_8$ |
| 2Ba) | $e$ | $b$ | $m$ | $a$ | $l$ | $p$ | $+1$ | $[\Pi/\Lambda] = (-4/pl)_8$ |

Note that, in case 1Ba), Lemma 13 gives $[\Pi/\Lambda] = (-4/l)_8(p/l)_4(l/p)_4$; but since $(p/l)_4(l/p)_4 = (-4/p)_8$ by Lemma 12, we get the relation in the table above.

As a matter of fact, the criteria involving $[\Pi/\Lambda]$ can just as well be obtained using genus theory (compare the discussion of $\mathcal{T}^{(\phi)}(2p)$ below). As the discussion of the $\phi$-part below shows, however, it seems that arguments from genus theory cannot always be replaced by the direct calculation of residue symbols.

$$\boxed{\text{The } \phi\text{-part}}$$

Our aim in this section is to show

**Proposition 14.** *If the torsor $\mathcal{T}^{(\phi)}(b_1)$ with $1 \neq b_1 \in \langle 2, p, l \rangle$ has a rational point, then the conditions in Table 2 must be satisfied.*

| $b_1$ | conditions (*) |
|---|---|
| 2 | $(-4/p)_8 = (-4/l)_8 = [\Pi/\Lambda] = 1$ |
| $p$ | $(p/l)_4 = (l/p)_4 = (-4/p)_8 = 1$ |
| $2p$ | $(p/l)_4(l/p)_4 = (-4/l)_8, (-4/p)_8 = 1, [\Pi/\Lambda] = (l/p)_4$ |
| $l$ | $(p/l)_4 = (l/p)_4 = (-4/l)_8 = 1$ |
| $2l$ | $(p/l)_4(l/p)_4 = (-4/p)_8, (-4/l)_8 = 1, [\Pi/\Lambda] = (p/l)_4$ |
| $pl$ | $(p/l)_4 = (l/p)_4, (-4/p)_8 = (-4/l)_8 = 1$ |
| $2pl$ | $(p/l)_4(l/p)_4 = (-4/p)_8 = (-4/l)_8, [\Pi/\Lambda] = 1$ |

TABLE 2. If $\mathcal{T}^{(\phi)}(b_1)$ has a rational point, then the conditions (*) must be satisfied.

For the proof of Prop. 14, we need the following proposition dealing with a slightly more general situation:

**Proposition 15.** *Let $k$ be a product of pairwise distinct primes $\equiv 1$ mod 8 that are quadratic residues of each other. Let $k = AB$ for $A, B \in \mathbb{N}$; if the torsor $\mathcal{T}^{(\phi)}(A)$ of $E_k$ has a nontrivial rational point, then the following conditions hold for any primary $\alpha \in \mathbb{Z}[i]$ with norm $A$:*

(1) $(-4/A)_8 = +1$;
(2) $[\alpha/\pi] = +1$ *for all* $\pi \mid B$;
(3) $(-4/p)_8 = (B/p)_4$ *for all* $p \mid A$.
(4) $[\alpha^*/\pi] = +1$ *for all* $\pi \mid \alpha$, *where* $\alpha = \alpha^*\pi$.

*Proof.* We have $\mathcal{T}^{(\phi)}(A) : AN^2 = M^4 + 4B^2e^4$; let $b = \gcd(M, B)$ be normalized by $b > 0$. Putting $N = bn$ and $M = bm$, we get $An^2 = b^2m^4 + 4c^2e^4$, where $bc = B$. We may assume that $m$ is odd: otherwise we switch the roles of $m$ and $e$. Note that $A \equiv 1$ mod 8 implies that $4 \mid e$.

Factoring the right hand side on $\mathbb{Z}[i]$ gives $\alpha\nu^2 = bm^2 + 2cie^2$ for some primary $\alpha \in \mathbb{Z}[i]$ with norm $N\alpha = A$. First observe that we have $\alpha\nu^2 \equiv bm^2 \equiv 1$ mod 8: thus $\alpha$ is congruent to a square modulo 8, and this implies 1. Moreover, $[\alpha/\pi] = [b/\pi] = (b/p) = +1$ for all $\pi \mid c$ with $N\pi = p$, and similarly $[\alpha/\pi] = 1$ for $\pi \mid b$, hence criterion 2.

Reducing the equation modulo some $\pi \mid \alpha$ gives $[1 + i/\pi](c/p)_4 = (-b/p)_4$, hence $(-4/p)_8 = (B/p)_4$ for all $p \mid A$, and this is 3.

Finally, subtracting $\alpha\nu^2 = bm^2 + 2cie^2$ from its conjugate yields $\alpha\nu^2 - \overline{\alpha}\overline{\nu}^2 = 4cie^2$; reducing modulo some $\overline{\pi} \mid \overline{\alpha}$ we get $[\alpha/\overline{\pi}] = (2c/p) = +1$. Since $[\pi/\overline{\pi}] = +1$, this is equivalent to $[\alpha^*/\pi] = +1$, proving 4. □

*Proof of Prop. 14.* In the case $\mathcal{T}^{(\phi)}(p)$ we have $A = p$ and $B = l$, so $(-4/p)_8 = 1$ from 1., $(p/l)_4(l/p)_4 = [\pi/\lambda] = 1$ from 2., $(-4/p)_8 = (l/p)_4$ from 3. and no condition from 4. In this way we find all criteria given in Table 2 except those involving $[\Pi/\Lambda]$. These have to be derived in an ad hoc manner:

• $\mathcal{T}^{(\phi)}(2)$ : $2n^2 = M^4 + p^2l^2e^4$. Write the torsor in the form $-p^2l^2e^4 = (M^2 + n\sqrt{2})(M^2 - n\sqrt{2})$. We assume that $(M, pl) = 1$; the other cases are treated similarly. Then $M^2 + n\sqrt{2} = \eta\Pi^2\Lambda^2\alpha^4$ for primes $\Pi, \Lambda \in \mathbb{Z}[\sqrt{2}]$ such that $N\Pi = p$, $N\Lambda = l$ and $\Pi \equiv \Lambda \equiv 1 \bmod 2$. Moreover, $\eta = \varepsilon^{\pm 1}$ with $\varepsilon = 1 + \sqrt{2}$. Adding the last equation to its conjugate gives $2M^2 = (\sqrt{2}M)^2 = \eta\Pi^2\Lambda^2\alpha^4 + \overline{\eta}\overline{\Pi}^2\overline{\Lambda}^2\overline{\alpha}^4$. Replacing $M$ by $M\varepsilon$ if necessary we may assume without loss of generality that $\eta = \varepsilon$. Thus

$$\varepsilon^{-1}(\sqrt{2}M)^2 \;=\; \Pi^2\Lambda^2\alpha^4 - \overline{\varepsilon}^2\overline{\Pi}^2\overline{\Lambda}^2\overline{\alpha}^4 \;=\; (\Pi\Lambda\alpha^2 + \overline{\varepsilon}\overline{\Pi}\overline{\Lambda}\overline{\alpha}^2)(\Pi\Lambda\alpha^2 - \overline{\varepsilon}\overline{\Pi}\overline{\Lambda}\overline{\alpha}^2).$$

Now $\Pi\Lambda\alpha^2 + \overline{\varepsilon}\overline{\Pi}\overline{\Lambda}\overline{\alpha}^2 \equiv \sqrt{2} \bmod 2$, hence $\Pi\Lambda\alpha^2 + \overline{\varepsilon}\overline{\Pi}\overline{\Lambda}\overline{\alpha}^2 = \sqrt{2}\mu^2$, $\Pi\Lambda\alpha^2 - \overline{\varepsilon}\overline{\Pi}\overline{\Lambda}\overline{\alpha}^2 = \sqrt{2}\varepsilon^{-1}\overline{\mu}^2$. Reducing modulo $\overline{\Pi}$ and using $[\Pi/\overline{\Pi}] = (2/p)_4$, $[\varepsilon/\overline{\Pi}] = (-4/p)_8 = 1$ (in this case), as well as $[\Lambda/\overline{\Pi}] = [\Lambda/\Pi]$ we find that the solvability of $\mathcal{T}^{(\phi)}(2)$ implies $[\Lambda/\Pi] = 1$.

• $\mathcal{T}^{(\phi)}(2p)$: Factoring $2pn^2 = M^4 + l^2e^4$ as $2pn^2 = (M^2 + le^2 + Me\sqrt{2l})(M^2 + le^2 - Me\sqrt{2l})$ and observing that $Me \equiv 1 \bmod 2$ implies that each factor is divisible exactly once by the prime ideal $\mathfrak{z}$ above 2. Thus $\mathfrak{z}\mathfrak{p}\mathfrak{n}^2 = (M^2 + le^2 + Me\sqrt{2l})$, where $\mathfrak{n}$ is an ideal with norm $n$. Let $h^+$ denote the class number of $\mathbb{Q}(\sqrt{2l})$ in the strict sense. We have to distinguish several cases:

(1) $h \equiv 2 \bmod 4$, $h^+ \equiv 4 \bmod 8$. By Proposition 2, this holds if and only if $(-4/l)_8 = -1$, and we also know that $N\varepsilon_{2l} = +1$ and that $\mathfrak{z}$ is principal in the wide sense. Now $[\Pi/\Lambda] = +1 \iff \mathfrak{p} \stackrel{+}{\sim} \boxed{4}$ by Lemma 3, and since $\mathfrak{z}\mathfrak{p}\mathfrak{n}^2$ is principal in the strict sense, this happens if and only if $\mathfrak{z}\mathfrak{n}^2 \stackrel{+}{\sim} \boxed{4}$. If $(2/l)_4 = -1$, then $\mathfrak{z} \stackrel{+}{\sim} 1$ is principal in the strict sense, and this happens if and only if $\mathfrak{n}^2 \stackrel{+}{\sim} \boxed{4}$, thus by genus theory $\iff (2/n) = (l/n) = +1$. But $(2/n) = (2p/l)_4 = -(p/l)_4$. Finally, solvability of $\mathcal{T}^{(\phi)}(2p)$ implies $(-4/l)_8 = (p/l)_4(l/p)_4$, so $(p/l)_4 = (-4/l)_8(l/p)_4 = -(l/p)_4$, and we see that $[\Pi/\Lambda] = (l/p)_4$ as claimed. If $(2/l)_4 = +1$, on the other hand, then $\mathfrak{z}$ is not principal in the strict sense, hence $[\Pi/\Lambda] = +1 \iff \mathfrak{n}^2 \stackrel{+}{\not\sim} \boxed{4}$, that is, iff $-1 = (2/n) = (p/l)_4$, and as above this gives $[\Pi/\Lambda] = (l/p)_4$.

(2) $h \equiv h^+ \equiv 4 \bmod 8$. By Proposition 2, this holds if and only if $(2/l)_4 = -1$ and $l \equiv 9 \bmod 16$. Here $\mathfrak{z}^2$ is principal in the strict sense and $\mathfrak{z}$ is not, in particular $\mathfrak{z} \stackrel{+}{\sim} \boxed{2}$ but $\mathfrak{z} \stackrel{+}{\not\sim} \boxed{4}$. Now $[\Pi/\Lambda] = +1 \iff \mathfrak{n} \stackrel{+}{\not\sim} \boxed{2}$ which in turn happens iff $-1 = (2/n) = (2p/l)_4 = -(p/l)_4$. Since $1 = (-4/l)_8 = (p/l)_4(l/p)_4$ from earlier solvability results, this gives $[\Pi/\Lambda] = 1 \iff (l/p)_4 = 1$ as claimed.

(3) $h^+ \equiv 0 \bmod 8$. By Proposition 2, this holds if and only if $(2/l)_4 = +1$ and $l \equiv 1 \bmod 16$. Here $\mathfrak{z}^2 = (2)$ is principal, and since the class group $\mathrm{Cl}_2^+(k)$ is cyclic, $\mathfrak{z} \stackrel{+}{\sim} \boxed{4}$. Thus $[\Pi/\Lambda] = +1 \iff \mathfrak{n} \stackrel{+}{\sim} \boxed{2} \iff 1 = (2p/l)_4 = (p/l)_4$, and we conclude as above that $[\Pi/\Lambda] = (l/p)_4$.

- $\mathcal{T}^{(\phi)}(2l) : N^2 = 2lM^4 + 2p^2le^4$. Symmetry reduces this to the discussion of $\mathcal{T}^{(\phi)}(2p)$.
- $\mathcal{T}^{(\phi)}(2pl) : N^2 = 2plM^4 + 2ple^4$. We start by factoring the torsor as $2pln^2 = M^4 + e^4 = (M^2 + e^2 + Me\sqrt{2})(M^2 + e^2 - Me\sqrt{2})$. Unique Factorization in $\mathbb{Z}[\sqrt{2}]$ gives $M^2 + e^2 + Me\sqrt{2} = \varepsilon\sqrt{2}\Pi\Lambda\nu^2$ and $M^2 + e^2 - Me\sqrt{2} = -\overline{\varepsilon}\sqrt{2}\,\overline{\Pi\Lambda}\,\overline{\nu}^2$. Subtracting the second equation from the first gives $2Me = \varepsilon\Pi\Lambda\nu^2 + \overline{\varepsilon}\,\overline{\Pi\Lambda}\,\overline{\nu}^2$, which in view of $[\varepsilon/\overline{\Pi}] = (-4/p)_8$ and $[\Pi/\overline{\Pi}] = (2/p)_4$ gives $[\Lambda/\Pi] = (Me/p)(-1/p)_8$.

On the other hand we have $2pln^2 = (M^2 + ie^2)(M^2 - ie^2)$, hence $M^2 + ie^2 = (1+i)\pi\lambda\nu^2$ for some $\nu \in \mathbb{Z}[i]$. This implies $(Me/p) = [Me/\pi] = [-i/\pi]_4 = (-1/p)_8$, hence our claim that $[\Pi/\Lambda] = 1$ is proved. □

The use of genus theory in this connection was suggested by the proofs of Pépin's conjectures in [20]. This concludes our discussion of the $\phi$-part of $\text{III}(\widehat{E}/\mathbb{Q})$.

## 5. The Main Result

The main result of this note is the following theorem:

**Theorem 16.** *Let $p \equiv l \equiv 1 \bmod 8$ be primes with $(p/l) = 1$. The properties of the Tate-Shafarevich groups $\text{III}(E_k/\mathbb{Q})[\phi]$ and $\text{III}(\widehat{E}_k/\mathbb{Q})[\psi]$ corresponding to the 2-isogenies between the elliptic curves $E_k : y^2 = x(x^2 - p^2l^2)$ and $\widehat{E}_k : y^2 = x(x^2 + 4p^2l^2)$ are recorded in Table 3. If the rank given there is 0, then the given subgroups actually equal $\text{III}(E_k/\mathbb{Q})[\phi]$ and $\text{III}(\widehat{E}_k/\mathbb{Q})[\psi]$, and we have $\text{III}(E/\mathbb{Q})[2] \simeq (\mathbb{Z}/2\mathbb{Z})^4$.*

The last two columns in Table 3 give the smallest examples of $p$ and $l$ satisfying the conditions and such that the given inequality for the rank is an equality (assuming the BSD-conjecture, and with the possible exception of the first line with $p = 41$, $l = 2273$, where the rank is 2 or 4). In all cases except one, the given example is the one that occurs first: the exception is $pl = 41 \cdot 1601$, where the example $pl = 41 \cdot 1321$ has the same residue symbols; yet rank $E_{41 \cdot 1321} = 0$.

Let us sketch the proof by going through one example. Take e.g. the second line; we claim that $\mathcal{T}^{(\phi)}(p)$ is the only possibly trivial torsor in $S^{(\phi)}(E/\mathbb{Q})$ (that means that it is the only one that might have a rational point). In fact, the torsors $\mathcal{T}^{(\phi)}(2)$, $\mathcal{T}^{(\phi)}(l)$, $\mathcal{T}^{(\phi)}(2l)$ and $\mathcal{T}^{(\phi)}(pl)$ are nontrivial since $(-4/l)_8 = -1$, whereas $\mathcal{T}^{(\phi)}(2p)$ and $\mathcal{T}^{(\phi)}(2pl)$ are nontrivial because $(p/l)_4(l/p)_4 \neq (-4/l)_8$. The other claims now follow immediately.

It remains to prove that $\text{III}(E/\mathbb{Q})[2]$ has order 16 if rank $E_{pl} = 0$. Recall the exact sequence

$$0 \longrightarrow \text{III}(E/\mathbb{Q})[\phi] \longrightarrow \text{III}(E/\mathbb{Q})[2] \longrightarrow \text{III}(\widehat{E}/\mathbb{Q})[\psi] \longrightarrow \widehat{C} \longrightarrow 0,$$

where $\widehat{C}$ is a finite 2-group of even rank by a result of Cassels. Since $\widehat{C}$ is a quotient of the group $\text{III}(\widehat{E}/\mathbb{Q})[\psi]$ of order 2 in our case, we must have $\widehat{C} = 0$, and in particular we get $\text{III}(E/\mathbb{Q})[2] \simeq \text{III}(E/\mathbb{Q})[\phi] \oplus \text{III}(\widehat{E}/\mathbb{Q})[\psi]$ as claimed.

The formula in [25, p. 30] shows that, assuming BSD, the order of $\text{III}(\widehat{E}/\mathbb{Q})[2]$ is 4 in these cases; this will be proved unconditionally in a subsequent paper dealing with a comparison of the method used here and classical 2-descent.

**Corollary 17.** *The curves of rank 0 among $E_{pl}$, where $p \equiv l \equiv 1 \bmod 8$ are primes such that $(p/l) = +1$, have density at least $\frac{1}{2}$. Those with rank 4 have density at most $\frac{1}{32}$.*

| $\left[\frac{\Pi}{\Lambda}\right]$ | $\left(\frac{l}{p}\right)_4$ | $\left(\frac{p}{l}\right)_4$ | $\left(\frac{-4}{p}\right)_8$ | $\left(\frac{-4}{l}\right)_8$ | $\mathbf{III}[\psi]$ | $\mathbf{III}[\phi]$ | rk E | $W^{(\phi)}$ | $p$ | $l$ |
|---|---|---|---|---|---|---|---|---|---|---|
| $+1$ | $+1$ | $+1$ | $+1$ | $+1$ | $1$ | $1$ | $\leq 4$ | $\langle 2, p, l\rangle$ | 41 | 2273 |
| | | | $+1$ | $-1$ | $1$ | $\langle 2p, l\rangle$ | $\leq 2$ | $\langle p\rangle$ | 41 | 769 |
| | | | $-1$ | $+1$ | $1$ | $\langle p, 2l\rangle$ | $\leq 2$ | $\langle l\rangle$ | 97 | 353 |
| | | | $-1$ | $-1$ | $\langle p\rangle$ | $\langle 2, p, l\rangle$ | $0$ | $1$ | 17 | 1361 |
| | | $-1$ | $+1$ | $+1$ | $1$ | $\langle p, l\rangle$ | $\leq 2$ | $\langle 2\rangle$ | 41 | 113 |
| | | | $+1$ | $-1$ | $1$ | $\langle p, l\rangle$ | $\leq 2$ | $\langle 2p\rangle$ | 113 | 233 |
| | | | $-1$ | $+1$ | $\langle p\rangle$ | $\langle 2, p, l\rangle$ | $0$ | $1$ | 17 | 953 |
| | | | $-1$ | $-1$ | $1$ | $\langle 2, p\rangle$ | $\leq 2$ | $\langle 2pl\rangle$ | 17 | 89 |
| | $-1$ | $+1$ | $+1$ | $+1$ | $1$ | $\langle p, l\rangle$ | $\leq 2$ | $\langle 2\rangle$ | 41 | 569 |
| | | | $+1$ | $-1$ | $\langle p\rangle$ | $\langle 2, p, l\rangle$ | $0$ | $1$ | 41 | 73 |
| | | | $-1$ | $+1$ | $1$ | $\langle p, l\rangle$ | $\leq 2$ | $\langle 2l\rangle$ | 17 | 457 |
| | | | $-1$ | $-1$ | $1$ | $\langle p, l\rangle$ | $\leq 2$ | $\langle 2pl\rangle$ | 17 | 433 |
| | | $-1$ | $+1$ | $+1$ | $\langle p\rangle$ | $\langle p\rangle$ | $\leq 2$ | $\langle 2, pl\rangle$ | 41 | 1601 |
| | | | $+1$ | $-1$ | $\langle p\rangle$ | $\langle 2, p, l\rangle$ | $0$ | $1$ | 41 | 449 |
| | | | $-1$ | $+1$ | $\langle p\rangle$ | $\langle 2, p, l\rangle$ | $0$ | $1$ | 17 | 569 |
| | | | $-1$ | $-1$ | $\langle p\rangle$ | $\langle 2, p, l\rangle$ | $0$ | $1$ | 17 | 977 |
| $-1$ | $+1$ | $+1$ | $+1$ | $+1$ | $\langle p\rangle$ | $\langle 2\rangle$ | $\leq 2$ | $\langle p, l\rangle$ | 113 | 569 |
| | | | $+1$ | $-1$ | $1$ | $\langle 2, l\rangle$ | $\leq 2$ | $\langle p\rangle$ | 41 | 433 |
| | | | $-1$ | $+1$ | $1$ | $\langle 2, p\rangle$ | $\leq 2$ | $\langle l\rangle$ | 17 | 353 |
| | | | $-1$ | $-1$ | $\langle p\rangle$ | $\langle 2, p, l\rangle$ | $0$ | $1$ | 73 | 89 |
| | | $-1$ | $+1$ | $+1$ | $\langle p\rangle$ | $\langle 2, p, l\rangle$ | $0$ | $1$ | 41 | 353 |
| | | | $+1$ | $-1$ | $\langle p\rangle$ | $\langle 2, p, l\rangle$ | $0$ | $1$ | 113 | 241 |
| | | | $-1$ | $+1$ | $1$ | $\langle 2, p\rangle$ | $\leq 2$ | $\langle 2l\rangle$ | 17 | 137 |
| | | | $-1$ | $-1$ | $\langle p\rangle$ | $\langle 2, p, l\rangle$ | $0$ | $1$ | 89 | 97 |
| | $-1$ | $+1$ | $+1$ | $+1$ | $\langle p\rangle$ | $\langle 2, p, l\rangle$ | $0$ | $1$ | 41 | 337 |
| | | | $+1$ | $-1$ | $1$ | $\langle 2, l\rangle$ | $\leq 2$ | $\langle 2p\rangle$ | 113 | 401 |
| | | | $-1$ | $+1$ | $\langle p\rangle$ | $\langle 2, p, l\rangle$ | $0$ | $1$ | 17 | 257 |
| | | | $-1$ | $-1$ | $\langle p\rangle$ | $\langle 2, p, l\rangle$ | $0$ | $1$ | 73 | 97 |
| | | $-1$ | $+1$ | $+1$ | $\langle p\rangle$ | $\langle p\rangle$ | $\leq 2$ | $\langle 2p, 2l\rangle$ | 113 | 257 |
| | | | $+1$ | $-1$ | $\langle p\rangle$ | $\langle 2, p, l\rangle$ | $0$ | $1$ | 41 | 241 |
| | | | $-1$ | $+1$ | $\langle p\rangle$ | $\langle 2, p, l\rangle$ | $0$ | $1$ | 89 | 257 |
| | | | $-1$ | $-1$ | $\langle p\rangle$ | $\langle 2, p, l\rangle$ | $0$ | $1$ | 17 | 281 |

TABLE 3. The Tate-Shafarevich groups $\mathbf{III}[\phi] := \mathbf{III}(E_k/\mathbb{Q})[\phi]$ and $\mathbf{III}[\psi] := \mathbf{III}(\widehat{E}_k/\mathbb{Q})[\psi]$ corresponding to the 2-isogenies between the elliptic curves $E_k : y^2 = x(x^2 - k^2)$ and $\widehat{E}_k : y^2 = x(x^2 + 4k^2)$ with $k = pl$ have subgroups as indicated. The column labeled by rk $E$ gives bounds for the rank of $E_k(\mathbb{Q})$. The column $W^{(\phi)}$ gives the subgroup of torsors in $S^{(\phi)}(E/\mathbb{Q})$ that may have rational points.

**Some Examples.** In [37], Wada and Taira (extending previous calculations of Noda & Wada [27]; see also Nemenzo [24]) computed the rank of most curves $E_k$ for $k < 40,000$. For 20 of these curves, they could only prove that the rank was between 2 and 4. For 8 out of these 20 numbers, our results show that the rank is in fact 2 in these cases:

| $k$ | $p$ | $l$ | $(l/p)_4$ | $(p/l)_4$ | $(-4/p)_8$ | $(-4/l)_8$ | $[\Pi/\Lambda]$ |
|---|---|---|---|---|---|---|---|
| 1513 | 17 | 89 | $+1$ | $-1$ | $-1$ | $-1$ | $+1$ |
| 2329 | 17 | 137 | $+1$ | $-1$ | $-1$ | $+1$ | $-1$ |
| 4633 | 41 | 113 | $+1$ | $-1$ | $+1$ | $+1$ | $+1$ |
| 6001 | 17 | 353 | $+1$ | $+1$ | $-1$ | $+1$ | $-1$ |
| 6953 | 17 | 409 | $+1$ | $+1$ | $-1$ | $+1$ | $-1$ |
| 7361 | 17 | 433 | $-1$ | $+1$ | $-1$ | $-1$ | $+1$ |
| 7769 | 17 | 457 | $-1$ | $+1$ | $-1$ | $+1$ | $+1$ |
| 9809 | 17 | 577 | $+1$ | $-1$ | $-1$ | $+1$ | $-1$ |

We remark in passing that the inequality rank $E \leq 2$ in these cases follows already from the criteria not involving $[\Pi/\Lambda]$. Moreover, the special case $k = 1513$ was discussed by Wada [36].

The tables of Nemenzo [25, 26] contain 70 more values $k = pl < 100,000$ such that $E_k$ has analytic rank 2 and Selmer rank 4. For 66 of them, the criteria involving the rational residue symbols suffice to show that the rank is at most 2; the 4 exceptions are $k = 64297 = 113 \cdot 569$, $67009 = 113 \cdot 593$, $93193 = 41 \cdot 2273$ and $94177 = 41 \cdot 2297$. For these values of $k$ we find $[\Lambda/\Pi] = -1$ except when $k = 93193$.

In [21], we will treat the remaining values of $k$ from [37] for which the rank could not be determined there.

CSU San Marcos, 333 S Twin Oaks Valley Rd, San Marcos, CA 92096-0001, USA
*E-mail address*: `franzl@csusm.edu`